\documentclass[12pt]{article}

\usepackage{amsmath,amsfonts,amsthm,amssymb, mathtools}
\usepackage{setspace}
\usepackage{ mathrsfs }
\usepackage{ marvosym }

\usepackage{Tabbing}
\usepackage{fancyhdr}
\usepackage{lastpage}
\usepackage{extramarks}

\usepackage{chngpage}
\usepackage{soul,color}
\usepackage{graphicx,float,wrapfig}
\usepackage{fullpage}
\usepackage{hyperref}
\usepackage[all]{xy}

\usepackage{enumerate}
\usepackage{marvosym}
\newcommand{\be}{\begin{enumerate}}
\newcommand{\ee}{\end{enumerate}}
\newcommand{\beq}{\begin{equation}}
\newcommand{\eeq}{\end{equation}}
\usepackage{verbatim}

\newcommand{\ba}{\begin{align*}}

\newcommand{\ea}{\end{align*}}

\newcommand{\R}{{\bf R}}

\newcommand{\Z}{{\bf Z}}

\newcommand{\RP}{{\mathbb RP}}

\newcommand{\ol}[1]{\overline{#1}}

\newcommand{\ra}{\rightarrow}
\newcommand{\tra}\twoheadrightarrow

\newcommand{\tla}\twoheadleftarrow

\newcommand{\vn}{\varepsilon}

\newcommand{\lp}{\left(}

\newcommand{\rp}{\right)}
\newcommand{\lpi}{\left|}
\newcommand{\rpi}{\right|}

\newcommand{\lbrac}{\left[}

\newcommand{\rbrac}{\right]}

\newcommand{\inj}{\operatorname{inj}}

\def\bea#1\eea{\begin{align*}#1\end{align*}}
\def\bc#1\ec{}

\newtheorem{Theorem}{Theorem}[section]

\newtheorem{Lemma}[Theorem]{Lemma}
\newtheorem{Proposition}[Theorem]{Proposition}
\newtheorem{Corollary}[Theorem]{Corollary}
\newtheorem{Conjecture}[Theorem]{Conjecture}

\let\inf\relax \DeclareMathOperator*\inf{\vphantom{p}inf}
\let\lim\relax \DeclareMathOperator*\lim{lim\vphantom{p}}

\title{Cyclic Pursuit on Compact Manifolds}
\date{}
\author{Dmitri Gekhtman}

\begin{document}

\maketitle

\abstract
\noindent We study a form of cyclic pursuit on Riemannian manifolds with positive injectivity radius. We conjecture that on a compact manifold,
the piecewise geodesic loop formed by connecting consecutive pursuit agents either collapses to a point in finite time or converges to a closed geodesic. The main result
is that this conjecture is valid for nonpositively curved compact manifolds.

\section{Introduction}
Our starting point is the classical three bug problem, first posed by Edouard Lucas \cite{L} in 1877: Three bugs start on the corners of an equilateral triangle, and each chases the next at unit speed. What happens? Answer: The bugs wind around the center of the triangle infinitely many times as they head inward along logarithmic spirals. They collide at the center of the triangle in finite time. We get similar behavior in general for a system of $n$ bugs starting at the vertices of a regular $n$-gon, each chasing its clockwise neighbor at unit speed (see e.g. \cite{Be}.) For illustrations of cyclic pursuit with initial conditions on a regular $n$-gon, see Clips 1 through 4 at the web address in the footnote.\footnote{\noindent \url{https://www.youtube.com/playlist?list=PLR5jDTSaPQj_xAkRpfipNLES_SSWIQeEp}} For details on the history of various versions of the $n$ bug problem, see the introduction of \cite{R1}.

Next, consider $n$ bugs starting at arbitrary positions in $\R^d$, with bug $i$ chasing bug $i+1$ mod $n$ at unit speed. Clips 5 through 7 at the address below demonstrate cyclic pursuit with randomly chosen initial conditions in the unit cube of $\R^3$. The typical observed behavior is as follows: 
Starting from the initial random configuration, chains of closely spaced bugs form, the chains come together to form a close approximation of a smooth knot, the knot unknots into an approximately circular loop, and the loop collapses to a point in finite time. The evolution of the piecewise linear loop connecting the bugs qualitatively
resembles the curve-shortening flow on the space of smooth loops in $\R^3$. In \cite{R1} and \cite{R2}, Richardson analyzed aspects of cyclic pursuit in $\R^d$. In \cite{R2}, he showed that for $n\geq 7$, the only stable configuration of $n$ bugs in cyclic pursuit is a planar regular $n$-gon. Based on simulations, Richardson conjectured that, if the initial positions of $n$ bugs are chosen uniformly at random in $[0,1]^d$, the probability of converging asymptotically to the regular $n$-gon configuration approaches 1 as $n\ra \infty$.

In this paper, we study cyclic pursuit on Riemannian manifolds with positive injectivity radius. To define pursuit in this case,
we choose initial positions such that each bug is within the injectivity radius of the next, and we have each bug chase the next with velocity equal to the unit vector tangent to the shortest geodesic connecting it to the next bug. Unlike in the Euclidean case, the bugs do not necessarily all
collide in finite time. Certainly, they cannot if the piecewise geodesic loop connecting consecutive
bugs is not null-homotopic. This leads to the conjecture that on a compact manifold, if the bugs
do not collide in finite time, the loop connecting them converges to a closed geodesic. The main result of this paper is that the conjecture is valid for pursuit on manifolds of nonpositive curvature. Clips 8 and 9 demonstrate convergence on a flat torus and flat M\"obius band, respectively. Numerical simulations suggest the conjecture is valid in general. Clip 9 demonstrates cyclic pursuit on $S^2$ and Clips 10, 11 show pursuit on $\RP^2$.     

The organization of the paper is as follows: In Section \ref{EUC}, we study basic properties of cyclic pursuit in Euclidean space. In Section \ref{Riemannian}, we introduce cyclic pursuit on Riemannian manifolds. In Section \ref{Convex}, we prove a result which states roughly that, if the bugs
enter a convex subset of a manifold, they stay in that subset. We derive as a consequence a condition for the pursuit to end in finite time. In Section \ref{Subsequence}, we prove subsequential convergence of the loop of bugs to a closed geodesic, and we obtain another criterion for pursuit to end in finite time. In Section \ref{Convergence}, we give
a condition for the loop of bugs to converge to a closed geodesic. In Section \ref{NEG}, we discuss convergence to closed geodesics which are locally length-minimizing, in the sense 
that any other loop uniformly close to the geodesic is longer. Then, we prove our {\em main result: for pursuit on a nonpositively curved compact manifold, the loop of bugs either collapses to a point in finite time or converges to a closed geodesic.}

\section{Notation}
\noindent
Unless otherwise stated, geodesics are parameterized at constant speed. 

If $(M, \langle\cdot, \cdot\rangle)$ is a Riemannian manifold and $p,q\in M$ are connected by a unique shortest geodesic, $[p,q]$ will denote the shortest geodesic from $p$ to $q$, parameterized as a map from the unit interval $[0,1]$. 

Unless otherwise stated, $\|\cdot\|:T_pM\ra \R$ will denote the norm associated to $\langle \cdot, \cdot \rangle_p$, and $d:M\times M\ra \R$ will denote the distance function associated to the metric. 

We identify $S^1\cong \R/\Z$. For each $x\in S^1$, we define $T_x:\R/\Z\ra \R/\Z$ to be the translation $T_x(s) = s+x$.

If $\alpha, \gamma: S^1\ra M$ are two loops, $d(\alpha,\gamma)$ will denote the supremum
distance $\sup_{s\in S^1}d(\alpha(s), \gamma(s))$. 

If $p\in M$ and $r>0$, $B_r(p)$ will denote the metric ball of radius $r$ centered at $p$.
For $p\in M$ and $K\subset M$, $d(p, K)$ will denote $\inf_{q\in K} d(p,q)$.
If $\alpha:S^1 \ra M$ is a closed geodesic and $\delta>0$, $\ol{N}_\delta(\alpha)$ will denote the closed $\delta$-neighborhood of $\alpha(S^1)$, $$\{p\in M | d(p, \alpha(S^1))\leq \delta\}.$$

\section{Cyclic Pursuit in Euclidean Space}\label{EUC}

We define cyclic pursuit of $n$ bugs in $\R^d$ with initial positions
$\{ b_i(0\}_{i\in \Z/n}$ as the unique collection of piecewise smooth functions
$\{b_i:[0,\infty) \ra \R^n\}_{i\in\Z/n}$ with the given initial conditions satisfying
\be
\item
If  $b_i(t) \neq b_{i+1}(t)$, then \beq \label{M1} \dot{b}_i(t) = \frac{b_{i+1}(t) - b_i(t)}{\|b_{i+1}(t) - b_i(t)\|}. \eeq
\item
If  $b_i(t_0) = b_{i+1}(t_0)$, then  $\forall t>t_0,~ b_i(t) = b_{i+1}(t)$.
\item
If  $b_i(t_0) = b_0(t_0)~\forall i$, then $\forall t>t_0, ~ b_i(t) = b_0(t_0).$
\ee

The following result is well known. We include a proof, as the proof will be useful later on.
\begin{Proposition}\label{FT}
For any set of initial conditions $\{b_i(0)\}_{i\in \Z/n}$, cyclic pursuit in $\R^d$ ends in finite time. That is, there is a $t_0>0$
so that $b_i(t) = b_0(t_0)$ for all $i$ and all $t\geq t_0$.
\end{Proposition}

\noindent {\em Proof:}

Let $l_i(t) = d(b_i(t), b_{i+1}(t))$, and let $l(t) = \sum_{i=1}^n l_i(t)$ be the length of the piecewise linear loop connecting the $b_i(t)$. 
We recall the following fact, which can be verified by direct computation:
Fix $p\in \R^d$ and consider the function $d_p:\R^d \setminus \{p\}\ra \R$,
$d_p(q) = d(p,q)$. Then the gradient of $d_p$ at $q$ is the unit vector $\frac{q-p}{\|q-p\|}$.

Let $u_i = \frac{b_{i+1}(t)-b_i(t)}{\| b_{i+1}(t)-b_i(t)\|}$.
Assuming for now that $l_i(t)>0$ for all $i$, we get $$\frac{d}{dt}l_i(t) = \langle u_i(t), \dot{b}_{i+1}(t)  \rangle + \langle -u_i(t), \dot{b}_i(t)\rangle  = 
\langle u_i(t), u_{i+1}(t)  \rangle + \langle -u_i(t), u_i(t) \rangle = \cos \theta_i(t)-1,$$
where $\theta_i(t) \in [0, \pi]$ is the angle between $u_i(t)$ and $u_{i+1}(t)$. By a theorem of Borsuk \cite{B}
, the sum of the exterior angles of a piecewise linear loop in $\R^d$ is at least $2\pi$.
So some $\theta_i$ is at least $\frac{2\pi}{n}$, and we find that $\frac{d}{dt}l(t) \leq \cos\frac{2\pi}{n}-1$.
In other words, $\frac{d}{dt}l(t)$ is negative, with absolute value bounded from below by 
$1-\cos\frac{2\pi}{n}$. If some $l_i(t)$ is 0, this effectively reduces $n$, so we still have the same bound on
$\frac{d}{dt}l(t)$. Thus, pursuit ends by time $l(0) \lbrac 1-\cos\frac{2\pi}{n} \rbrac^{-1}$.\qed

{\em Remark}~:
If the $\theta_i$ are all $\leq \frac{\pi}{2}$, Jensen's inequality applied to
$1-\cos \theta$ on $[0, \frac{\pi}{2}]$ yields
$\lpi  \frac{d}{dt} l(t) \rpi \geq n(1- \cos\frac{2\pi}{n})$.
On the other hand, if at least one of the $\theta_i$ is greater than $\frac{\pi}{2}$,
then $\lpi  \frac{d}{dt} l(t) \rpi\geq 1-\cos\frac{\pi}{2} = 1$.
Thus, assuming $l_i(t)>0$ for all $i$, we have $\lpi  \frac{d}{dt} l(t) \rpi\geq  \min\lbrac 1, n(1- \cos\frac{2\pi}{n})\rbrac$.
Since $\min\lbrac 1, n(1- \cos\frac{2\pi}{n})\rbrac$ is a nonincreasing function of $n\geq 2$, we still have $\lpi  \frac{d}{dt} l(t) \rpi\geq \min\lbrac 1, n(1- \cos\frac{2\pi}{n})\rbrac$
if some (but not all) of the $l_i(t)$ are 0.
Hence, the time from the start of the pursuit process to its end
is bounded above by $l(0)\cdot \lp \min\lbrac 1, n(1- \cos\frac{2\pi}{n})\rbrac\rp^{-1}$, which grows linearly in $n$.
(Compare this to the $O(n^2)$ bound on the time obtained from the estimate 
$\lpi \frac{d}{dt}l(t)\rpi \geq 1-\cos\frac{2\pi}{n}$ in the last paragraph.
)
Note also that, in the case that the $b_i(0)$ are vertices of a regular planar $n$-gon,
we get that the time to mutual capture is precisely $l(0)\cdot\lbrac n(1-\cos\frac{2\pi}{n}) \rbrac^{-1}$. 

\section{Pursuit on Riemannian Manifolds}\label{Riemannian}
Cyclic pursuit on a Riemannian manifold is defined just as in the
Euclidean case: each bug's velocity is the unit vector pointing towards
the next bug along the shortest geodesic connecting the two. To ensure
that there is a unique shortest geodesic connnecting each pair of
bugs, we consider only manifolds with positive injectivity
radius, and we choose initial positions so that the distance between each
bug and its prey is less than the injectivity radius.

Let $(M,g)$ be a manifold with positive injectivity radius $\text{inj}\left(M\right)$,
and let $\left\{ b_{i}(0)\right\} _{i\in\mathbb{Z}/n}$ be initial positions
in $M$ satisfying $d\left(b_{i}(0),b_{i+1}(0)\right)<\text{\text{inj}}\left(M\right).$
Then we define $ \{ b_{i}:[0,\infty)\rightarrow M \} _{i\in\mathbb{Z}/{n}}$
as the unique collection of piecewise smooth functions with the given initial conditions
satisfying
\be
\item
If  $b_i(t) \neq b_{i+1}(t)$, then  \begin{equation} \label{M2} \dot{b}_i(t) = \frac{ \exp_{b_i(t)}^{-1} (b_{i+1}(t) ) }{\| \exp_{b_i(t)}^{-1} (b_{i+1}(t) )\|}.\end{equation}
\item
If  $b_i(t_0) = b_{i+1}(t_0)$, then  $\forall t>t_0,~ b_i(t) = b_{i+1}(t)$.
\item
If  $b_i(t_0) = b_0(t_0)~\forall i$, then $\forall t>t_0, ~ b_i(t) = b_0(t_0).$
\ee
Let $l_i(t) = d(b_i(t), b_{i+1}(t))$.
To see that the pursuit process is well-defined for all $t\geq 0$, we need to check that each
$l_i(t)$ is non-increasing and thus stays less than $\inj(M)$.

To compute $\frac{d}{dt}l_i(t)$, we recall the following fact, which follows from the Gauss lemma of Riemannian geometry:
If $p\in M$ and $U$ is a normal neighborhood of $p$, consider the function $d_p:U\setminus \{p\}\ra \R$ given by $d_p(q) = d(p,q)$.
The gradient of $d_p$ at $q$ is the tangent at $q$ of the shortest unit speed geodesic going from $p$ to $q$.

Now, if $l_i(t)>0$, we use the above fact and the law of motion \eqref{M2} to compute, just as in the last section, that
$$\frac{d}{dt}l_i(t) = \cos \theta_i(t)-1,$$
where $\theta_i$ is the angle at $b_{i+1}(t)$ between $[b_{i}(t),b_{i+1}(t)]$ and $[b_{i+1}(t),b_{i+2}(t)]$. 
So if $l_i(t)>0$, then $\frac{d}{dt}l_i(t)\leq 0$ and thus $l_i$ is locally non-increasing at $t$.
On the other hand, if $l_i(t)=0$, then $l_i(t')=0$ for all $t'>t$. So each $l_i$ is indeed non-increasing,
and the pursuit process is well-defined.

For each $t\geq 0$, $i\in \Z/n$, let $\beta^t_i = [b_i(t), b_{i+1}(t)]$ be the shortest geodesic connecting $b_i(t)$ to $b_{i+1}(t)$.
Let $\beta^t:\R/\Z \ra M$ be the constant-speed piecewise geodesic loop formed by concatenating
the $\beta^t_i$, with $\beta^t(0) = \beta^t(1) = b_0(t)$. Then $t\mapsto \beta^t$ is a homotopy of
loops, so if $\beta^0$ is not null-homotopic, the pursuit process will not end in finite time.
So Proposition \ref{FT} does not generalize to pursuit on Riemannian manifolds.
This leads to the following conjecture for compact manifolds:
\begin{Conjecture}\label{C}
If $M$ is a compact Riemannian manifold, and $\{b_i(0)\}_{i\in \Z/n}$ are initial conditions for cyclic pursuit,
then the associated family of loops $\beta^t$ either collapses to a point in finite time or converges to a closed geodesic as $t\ra \infty$.
\end{Conjecture}
By convergence above, we mean convergence in the quotient of $C^0(S^1, M)$ by rotations in the domain. In other words, a sequence of loops $\{  \gamma_j \}_{j=1}^\infty$ converges to $\gamma:\R/\Z\ra M$
if 
\begin{equation}\label{CONV}
 \lim_{j\ra \infty} \adjustlimits\inf_{c\in \R} \sup_{s\in \R/\Z} d(\gamma_j(s), \gamma(s+c)) =0.
\end{equation}

We prove Conjecture \ref{C} in the case of pursuit on nonpositively curved compact manifolds in Section \ref{NEG}.

{\em Remark}~: As observed above, if pursuit ends in finite time, then $\beta^0$ is nullhomotopic. The converse is not true.
For instance, suppose $\alpha$ is a nullhomotopic closed geodesic along which all sectional curvatures are negative (e.g. the neck of a dumbbell.)
We will see in Section \ref{NEG} that if $\beta^0$ is sufficiently close to $\alpha$, then $\beta^t$ will converge to $\alpha$.



\section{Convex Submanifolds}\label{Convex}
We will need the following result, which states roughly that, if at some time the $b_i$ all belong to a convex set $K\subset M$,
then they stay in $K$.
\begin{Proposition}\label{CO}
Let $M^d$ be a Riemannian manifold with $\inj(M)>0$,
$\{b_i(t)\}_{i\in \Z/n}$ cyclic pursuit curves on $M$, $l_i(t) = d(b_i(t), b_{i+1}(t))$.
Let $K^d\subset M$ be a smoothly embedded submanifold with boundary, topologically closed in $M$. Suppose there is an $R \in \lp 0,\inj(M)\rp$ so that for any two points $p_1,p_2 \in K$ with $d(p_1,p_2)<R$, the geodesic segment $[p_1,p_2]$ is contained in $K$.  
If for some $t_0\in [0,\infty)$, all of the $b_i(t_0)$ are in $K$ and all of the $l_i(t_0)$ are less than $R$, then
$b_i(t)\in K$ for all $i\in \Z/n$, $t\geq t_0$.  
\end{Proposition}

\noindent {\em Proof idea:}
If one of the bugs reaches $\partial K$, then, by the convexity assumption, the bug's velocity will not point out of $K$. So the bug will stay in $K$. 

\smallskip

\noindent {\em Proof:}
Since $K$ is closed, embedded, and of the same dimension as $M$, its topological boundary in $M$ is the boundary manifold $\partial K$.
Suppose for the sake of contradiction that there is a $t_1>t_0$ and $j\in \Z/n$ so that $b_j(t_1)$ is not in $K$.
Set $$t' = \sup \{ t\in [t_0, t_1] |b_i(t)\in K~\forall i \},~ \vn = t_1-t'.$$
Then since $K$ is closed and the $b_i$ are continuous, all of the $b_i(t')$ are in $K$ and at least one of the $b_i(t')$ is in $\partial K$.
Furthermore, for all $t\in (t', t'+\vn]$, at least one of the $b_i(t)$ is in $M\setminus K$.

For each $i\in \Z/n$, let $(V_i, x_i^1,\ldots, x_i^d)$ be a coordinate neighborhood of $b_i(t')$
with the property that $V_i\cap K = \{p\in V_i| x_i^d(p)\leq 0\}$, and thus that $V_i\cap \partial K = \{p\in V_i| x_i^d(p)= 0\}$.
(If $b_i(t')$ is in the interior of $K$, it may be that $V_i\cap \partial K = \emptyset$ and $x_i^d<0$ on all of $V_i$.)
Shrinking $\vn$ if necessary, we may assume $b_i([t',t'+\vn])\subset V_i$ for all $i$.
Let $b_i^d=x^d_i \circ b_i$ denote the $d$-th component of $b_i$.
Let $h(t) = \max_i b_i^d(t)$.
Since all of the $b_i(t')$ are in $K$, and at least one of $b_i(t')$ is in $\partial K$, $h(t')=0$. 
For each $t\in (t', t'+\vn]$, at least one of the $b_i(t)$ is in the complement of $K$, so $h(t)>0$.
Assume without loss of generality that $b_i(t')\neq b_{i+1}(t')$, for all $i$. Then take $\vn$ small enough that $b_i(t)\neq b_{i+1}(t)$, 
for all $t\in [t', t'+\vn]$ and all $i\in \Z/n$. Since each $b_i^d$ is smooth on $[t', t'+\vn]$, $h$ is absolutely continuous on $[t', t'+\vn]$.
So $h$ is almost everywhere differentiable, and we have $h(t) = \int_{t'}^t \frac{d}{ds} h(s) ds$ for each $t\in [t', t'+\vn]$.
Thus, for some $c_1 \in [0,\vn]$, we have that $\frac{d}{dt}h(t'+c_1)$ is defined and
$$h(t'+\vn) \leq \vn \frac{d}{dt}h(t'+c_1)= \vn \frac{d}{dt}b_j^d(t'+c_1),$$
for some $j$ for which $b_j^d(t'+c_1) = h(t'+c_1)$. Take $\vn$ small enough that if $b_i^d(t) = h(t)$
for some $t\in [t', t'+\vn]$, $b_i^d(t') = 0$. Then in particular, $b^d_j(t')=0$ for the $j$ in the last displayed formula.

For each $i\in \Z/n$, let $v_i^d(p,q)$ be the $d$-th component of the initial unit tangent to $[p, q]$, for $(p,q)\in V_i\times V_{i+1}$ with $0<d(p,q)<\inj(M)$.  
By the law of motion \eqref{M2},
$$\frac{d}{dt}b_j^d(t'+c_1)=v^d_j(b_j(t'+c_1), b_{j+1}(t'+c_1)). $$
For each $i$, let $B_i\subset V_i$ be an open coordinate ball centered at $b_i(t')$ with $\ol{B}_i\subset V_i$.
Shrinking the $B_i$ if necessary, assume there is a $\delta>0$ so that $\delta<d(p,q)<R$ for all $(p,q)\in B_i\times B_{i+1}$. Then $v^d$
is $C^1$ on $\ol{B}_i\times \ol{B}_{i+1}$. Since $K$ contains $[p,q]$ whenever $p,q\in K$ and $d(p,q)<R$,
we have for all $(p,q)\in B_i\times B_{i+1}$ with $x_i^d(p)=0$ and $x_{i+1}^d(q)\leq 0$, that $v_i^d(p,q)\leq 0$.
Let $\mu_i$ (resp. $\nu_i$) be the maximum on $\ol{B}_i\times \ol{B}_{i+1}$ of the absolute value of the derivative of $v_i^d(p,q)$ with respect to the $d$-th component
of $p$ (resp. $q$.) Let $\mu = \max_i \mu_i$, $\nu = \max_i = \nu_i$. 
Taking $\vn$ small enough that $b_{i}\lp [t',t'+\vn]\rp\subset B_i$ for all $i$,
we have
$$v^d_j(b_j(t'+c_1), b_{j+1}(t'+c_1)) \leq (\mu+\nu) h(t'+c_1).$$

From the last three displayed formulas, we get
$$h(t'+\vn) \leq (\mu+\nu) \vn h(t'+c_1).$$
Similarly, $h(t'+c_1) \leq (\mu+\nu)\vn h(t'+c_2)$
for some $c_2\in [0,c_1]$, so we obtain
$h(t'+\vn) \leq \lp (\mu+\nu)\vn\rp^2  h(t'+c_2)$.
Inductively, we get for each positive integer $k$,
$$h(t'+\vn) \leq \lp (\mu+\nu)\vn\rp^k   h(t'+c_k)$$
for some $c_k \in [0,\vn]$. Let $C=\max_{i\in \Z/n}\sup_{p\in B_i} \lpi x_i^d(p)\rpi$. Taking $\vn<\frac{1}{2}(\mu+\nu)^{-1}$,
we get $h(t'+\vn) < 2^{-k}C$. Letting $k\ra \infty$ yields $h(t'+\vn)= 0$, a contradiction.

\qed

We say a subset $K$ of a Riemannian manifold $M$ is {\em convex}, if for each pair $p,q\in K$, there is a unique shortest geodesic in $M$ connecting $p,q$
and this geodesic is contained in $K$. Proposition \ref{CO} is the key ingredient in the proof of the following result:
\begin{Proposition}\label{FT2}
Suppose $M$ is compact, and $\{b_i(t)\}_{i\in \Z/n}$ are cyclic pursuit curves in $M$. Let $l_i(t) = d(b_i(t), b_{i+1}(t))$. If $l_i(t)\ra 0$ for all $i$, then pursuit ends in finite time. 
\end{Proposition}

\noindent {\em Proof idea:}
Reduce to the Euclidean case by noting that the bugs will eventually lie in a small, convex, approximately Euclidean ball. 

\noindent {\em Proof:}

Assume for the sake of contradiction that pursuit does not end in finite time. Without loss of generality, assume $l_i(t)>0$ for all $i$.
Since $M$ is compact, there is a $p\in M$ and a sequence of times $(t_j)_{j=1}^\infty$ with $t_j\ra \infty$ so that $b_0(t_j)\ra p$ as $j\ra \infty$.
Let $r\in (0,\inj(M))$ be small enough that the closed $r$-ball centered at $p$, $\ol{B}_r(p)$, is convex.
Since $d(b_i(t_j), b_{i+1}(t_j))\ra 0$ and $b_0(t_j)\ra p$,
$b_i(t_j) \ra p$ for all $i$, so there is an $J$ for which all of the $b_i(t_J)$ belong to $\ol{B}_r(p)$.
By Proposition \ref{CO}, the $b_i(t)$ remain in $\ol{B}_r(p)$ for all $t>t_J$. 

Let $(U, x^i)$ be a normal coordinate neighborhood centered at $p$. By Corollary \ref{AN} in the Appendix, we have that for small enough
$r$, $\ol{B}_r(p)$ is a convex subset of $U$ and has the following property: for any two geodesics
$\gamma_1:[0,a_1] \ra \ol{B}_r(p),\gamma_2:[0,a_1] \ra \ol{B}_r(p)$ with $\gamma_1(0) = \gamma_2(0)$, the metric angle between 
$\dot{\gamma}_1(0)$ and $\dot{\gamma}_2(0)$ is within $\frac{\pi}{n}$ of the Euclidean angle,
computed in the coordinates $(U,x^i)$, between $\gamma_1(a_1)-\gamma_1(0)$ and $\gamma_2(a_2)-\gamma_2(0)$. 

Now, choose $r$ as above
and find $t_0$ so that all of the $b_i(t)$ are in $\ol{B}_r(t)$ for $t\geq t_0$. 
As we showed in the proof of Proposition \ref{FT}, at least one of the Euclidean angles
of the piecewise linear loop connecting the $b_i(t)$ is $\geq \frac{2\pi}{n}$.
So by the result quoted the last paragraph, at least one of the angles of the piecewise geodesic loop 
connecting the $b_i(t)$ is $\geq \frac{\pi}{n}$. Thus, 
$ \frac{d}{dt}l(t) \leq \cos\frac{\pi}{n}-1$ for $t\geq t_0$
and so pursuit ends by time $t_0 + l(t_0)(1- \cos\frac{\pi}{n})^{-1}$. This is a contradiction.
\qed

\section{Subsequential Convergence}\label{Subsequence}
In this section, $M$ is a compact Riemannian manifold.

\begin{Proposition}\label{SS}
Let $\{b_i(t)\}_{i\in \Z/n}$ be pursuit curves on $M$, $\beta^t$ be the associated family of piecewise geodesic loops, $l(t)$ the length of $\beta^t$. If the pursuit does not end in finite time, then there is a sequence of times $(t_j)_{j=1}^\infty$, $t_j\ra \infty$
so that $\beta^{t_j}$ converges uniformly to a closed geodesic of length $L=\lim_{t\ra\infty}l(t)$ as $j\ra \infty$.
\end{Proposition}

\noindent {\em Proof sketch:}
Take a sequence $t_j$ so that the $b_i(t_j)$ converge and so that $\frac{d}{dt}l(t_j)\ra 0$. Then $\beta^{t_j}$ converges to
a piecewise geodesic loop. The condition $\frac{d}{dt}l(t_j)\ra 0$ implies that the angles between segments of the limiting loop are 0.

\noindent {\em Proof:}

Let $\beta^t_i = [b_i(t), b_{i+1}(t)]$. Recall that $\beta^t$ is the constant-speed piecewise geodesic loop formed from the $\beta^t_i$, with $\beta^t(0) = b_0(t)$. 
We have $L>0$, else by Proposition \ref{FT2}, pursuit ends in finite time.
Assume without loss of generality that $d(b_i(t), b_{i+1}(t))>0$ for all $i\in \Z/n,t>0$. Then we have for all $t$ that
$$\frac{d}{dt}l(t) = \sum_i \lp \cos\theta_i(t)-1\rp,$$ where $\theta_i(t)$ is the angle at $b_{i+1}(t)$
between $\beta^t_i$ and $\beta^t_{i+1}$. Since $l$ is differentiable, nonincreasing, and bounded from below,
there is a sequence $(t_j)_{j=1}^\infty$, $t_j\ra \infty$, so that $\frac{d}{dt}l(t_j) \ra 0$.
This implies that for each $i$,  $\theta_i(t_j)\ra 0$ as $j\ra \infty$. Since $M$ is compact, we may pass to a subsequence
and assume that for each $i$, $b_i(t_j)$ converges to some point $a_i \in M$. Then $[b_i(t_j), b_{i+1}(t_j)]$ converges uniformly to $[a_i, a_{i+1}]$.
Let $\alpha$ be the constant speed piecewise geodesic loop
formed from the geodesic segments $[a_i, a_{i+1}]$, with $\alpha(0) = a_0$. Then $\beta^{t_j}$ converges to $\alpha$
uniformly. By continuity,
$\alpha$ has length $L$.

We need to show that $\alpha$ is a closed geodesic. To do this, it suffices to show that the angles between successive geodesic
segments comprising $\alpha$ are 0. We need to include the case that $a_i = a_{i+1}$ for some $i$.
To this end, suppose $a_{i-1} \neq a_i$, and let $k$ be the smallest integer so that
$a_i = a_{i+1} = \cdots = a_{i+k}$. We need to show that the angle at $a_i$ between $[a_{i-1}, a_i]$ and $[a_i, a_{i+k+1}]$ is 0.
Let $(U, x^i)$ be a normal coordinate neighborhood centered at $a_i$, and $\|\cdot\|_U$ be the Euclidean norm on $TU$ coming from the coordinates. Fix $\vn>0$. Then for large enough $j$, $b_i(t_j),\ldots, b_{i+k}(t_j)$
are in $U$ and $$\left \| \frac{\dot{\beta}^{t_j}_m(1)}{\|\dot{\beta}^{t_j}_m(1) \|} - \frac{\dot{\beta}^{t_j}_m(0)}{\|\dot{\beta}^{t_j}_m(0) \|} \right\|_U<\vn$$
for $m=i,\ldots, i+k-1$. (See formula \eqref{DIR} of the Appendix.) Since $\theta_i(t_j)\ra 0$ for all $i$, we have for large enough $j$ that 
$$\left \| \frac{\dot{\beta}^{t_j}_{m}(0)}{\|\dot{\beta}^{t_j}_{m}(0) \|} - \frac{\dot{\beta}^{t_j}_{m-1}(1)}{\|\dot{\beta}^{t_j}_{m-1}(1) \|} \right\|_U<\vn$$
for $m=i,\ldots, i+k$. From the last two displayed expressions, we obtain
$$\left \| \frac{\dot{\beta}^{t_j}_{i+k}(0)}{\|\dot{\beta}^{t_j}_{i+k}(0) \|} - \frac{\dot{\beta}^{t_j}_{i-1}(1)}{\|\dot{\beta}^{t_j}_{i-1}(1) \|} \right\|_U<(2k-1)\vn.$$
Thus, the expression on the left of the last inequality converges to 0 as $j\ra \infty$.
But  $ \frac{\dot{\beta}^{t_j}_{i+k}(0)}{\|\dot{\beta}^{t_j}_{i+k}(0)\|}$ converges to the unit tangent to $[a_i, a_{i+k+1}]$ at $a_i$,
and   $\frac{\dot{\beta}^{t_j}_{i-1}(1)}{\|\dot{\beta}^{t_j}_{i-1}(1)\|}$ converges to the unit tangent to $[a_{i-1}, a_i]$ at $a_i$.
Hence, these two unit tangent directions are the same, and the angle between $[a_i, a_{i+k+1}]$ and $[a_{i-1}, a_i]$ at $a_i$ is 0,
as claimed. \qed

As a consequence of the last Proposition, we have
\begin{Corollary}
If for some $t_0 \geq 0$, the length of $\beta^{t_0}$ is less than the length $\lambda_{\min}$ of the shortest closed geodesic of $M$, pursuit ends in finite time.
In particular, if the length of $\beta^0$ is less than $\lambda_{\min}$, pursuit ends in finite time.
\end{Corollary} \qed

We also have
\begin{Corollary}
If for some $t_0 \geq 0$, the $b_i(t_0)$ all lie in a convex, smoothly embedded, closed metric ball $\ol{B}\subset M$, then pursuit ends in finite time.
\end{Corollary}

\noindent {\em Proof:}

By Proposition \ref{CO}, the $b_i(t)$ stay in $\ol{B}$ for $t>t_0$. If pursuit does not end in finite time, then arguing as in Proposition \ref{SS}, there is a sequence
$t_j\ra \infty$ so that $\beta^{t_j}$ converges to a closed geodesic in $\ol{B}$. But there are no closed geodesics contained in $\ol{B}$. 

\qed

It follows, for example, that if the $b_i(t_0)$ all lie in an open hemisphere of $S^2$ with its standard metric, pursuit ends in finite time. 

\section{A Criterion for Convergence}\label{Convergence}

The next result gives a criterion for convergence of $\beta^t$ to a closed geodesic $\alpha$. 
\begin{Proposition}\label{DONE}
Let $M$ be a Riemannian manifold with $\inj(M)>0$. Let $\{b_i(t)\}_{i\in \Z/n}$ be cyclic pursuit curves on $M$, $\beta^t$ the associated family of loops.
Suppose there is a sequence $t_j\ra \infty$ and a closed geodesic $\alpha$ so that $\beta^{t_j}\ra \alpha$ uniformly.
 If $\sup_{s\in S^1}d(\beta^t(s), \alpha(S^1))\ra 0$ as $t \ra \infty$, then
$\beta^t$ converges to $\alpha$ in the sense of Equation $\ref{CONV}$, as $t\ra \infty$.   
\end{Proposition}

\noindent {\em Proof idea:}
For large $t$, the curves $\beta^t$ and $\alpha$ have approximately the same length. In addition, $\beta^t$ fits into a small tubular neighborhood of $\alpha$.
These two facts force $\beta^t$ to be uniformly close to $\alpha$. 

\noindent {\em Proof:}

Let $U$ be an open neighborhood of $\alpha(S^1)$ such that each $p\in U$ has a unique closest point $\pi(p)$ on $\alpha(S^1)$. 
Shrinking $U$ if necessary, we may construct a smooth unit vector field $X$ on $U$ extending the unit tangent field 
$\frac{\dot{\alpha}}{\|\dot{\alpha}\|}$ of $\alpha$.

Fix $\vn>0$. Let $l_i(t) = d(b_i(t), b_{i+1}(t))$, and let $\lambda$ be the minimum of $\lim_{t\ra\infty} l_i(t)$ over $i$ for which $\lim_{t\ra\infty} l_i(t)>0.$ The following fact
follows from the continuous dependence of the initial unit tangent of a geodesic $[p,q]$ on the endpoints $p$ and $q$:
There is a $\delta$ so that $\ol{N}_\delta(\alpha)\subset U$ and if $\gamma:[0,1] \ra \ol{N}_\delta(\alpha)$ is a geodesic
of length $\geq \lambda$, then the component of $\frac{\dot{\gamma}(0)}{\|\dot{\gamma}(0)\|}$ normal to
$X(\gamma(0))$ has length less than $\vn$.

Consider $i$ such that  $\lim_{t\ra\infty} l_i(t)>0$. Let $a^t_i = \pi( \beta^t_i(0) )$. By hypothesis, 
\begin{equation}\label{Y}
d\lp \beta^t_i(0), a^t_i  ) \rp \ra 0
\end{equation}
as $t\ra \infty$.
By the observation in the previous paragraph, the component of $\dot{\beta}^t_i(0)$ orthogonal to $X(\beta^t_i(0) )$ goes to 0 as well.
Since $\beta^{t_j}\ra \alpha$ uniformly, we have by continuity that 
\begin{equation}\label{Z}
\lpi\lpi \dot{\beta}^t_i(0) - l_i(t)\cdot X(\beta^t_i(0))\rpi\rpi \ra 0.
\end{equation}

Now, let $\alpha_i^t$ be the segment of $\alpha$ starting at $a_i^t$ with initial velocity $ l_i(t)\cdot X(a_i(t))$.
Since a geodesic depends continuously on its initial parameters, \eqref{Y} and \eqref{Z} give 
\begin{equation}\label{X}
\lim_{t\ra \infty}\sup_{s\in [0,1]} d\lp \beta^t_i(s), \alpha_i^t(s)\rp= 0.
\end{equation} 

Let $i_1,\ldots, i_m$ be the values of $i$ for which $\lim_{t\ra \infty} l_i(t)>0$, listed in order.
Now, let $\gamma^t:\R/\Z\ra M$ be the piecewise continuous loop formed by concatenating the segments $\alpha^t_{i_j}\lp  [0,1) \rp$.
We parameterize $\gamma^t$ so that each $\alpha^t_{i_j}$ is traversed at the same constant speed, and $\gamma^t(0) = \alpha^t_{i_1}(0)$.
As a consequence of \eqref{X}, 
\begin{equation}\label{U} 
\lim_{t\ra \infty}\sup_{s\in S^1} d\lp \beta^t(s), \gamma^t(s) \rp = 0.
\end{equation}

Let $c^t$ be such that $\alpha(c^t) = \gamma^t(0)$. By triangle inequality, $\lim_{t\ra\infty} d( \alpha^t_{i_j}\lp  1 \rp, \alpha^t_{i_{j+1}}\lp  0 \rp)= 0$ for $j=1,\ldots, m-1$.
Also, $\sum_{j=1}^m l_{i_j}(t)$ converges to the length of $\alpha$ as $t\ra \infty$. It follows that
\begin{equation}\label{V} 
\lim_{t\ra \infty}\sup_{s\in S^1} d\lp \gamma^t(s), \alpha(s+c^t) \rp = 0.
\end{equation}

From \eqref{U} and \eqref{V}, we get
$$
\lim_{t\ra \infty}\sup_{s\in S^1} d\lp \beta^t(s), \alpha(s+c^t)\rp = 0,
$$
which completes the proof.

\section{Nonpositive Curvature}\label{NEG} 
In the next proposition, we show that if a subsequence $\beta^{t_j}$ converges to a closed geodesic $\alpha$ which is
a local minimizer of length, then $\beta^t$ converges to $\alpha$.

We recall the following notation: For each $x\in S^1$, let $T_x:\R/\Z\ra \R/\Z$ be the translation $T_x(s) = s+x$.
If $\alpha, \gamma: S^1\ra M$ are two loops, $d(\alpha,\gamma)$ denotes the supremum
distance $\sup_{s\in S^1}d(\alpha(s), \gamma(s))$. 

\begin{Proposition}\label{ISO}
Let $M$ be a Riemannian manifold with $\inj(M)>0$, $\{b_i(t)\}_{i\in \Z/n}$ cyclic pursuit curves on $M$, $\beta^t$ the associated family of loops.
Let $\alpha$ be a closed geodesic of length $L$, and suppose there is $\vn>0$ so that any rectifiable, constant speed loop $\gamma$ of length $L$
with $d(\alpha,\gamma)< \vn$ is a reparameterization of $\alpha$, i.e. $\gamma = \alpha \circ T_x$ for some $x\in S^1$.  
If there is a sequence of times $t_j\ra\infty$ so that $\beta^{t_j}\ra \alpha$ as $j\ra \infty$, then $\beta^t\ra \alpha$ as $t\ra \infty$.
\end{Proposition}

\noindent{\em Proof sketch:}
If $\beta^t$ does not converge to $\alpha$, then $\beta^t$ has another length $L$ subsequential limit $\gamma$.
It follows that for any $k$, there is a homotopy from $\alpha$ to $\gamma$ through curves of length between
$L$ and $L+\frac{1}{k}$. The homotopy passes through a curve $\eta_{k}$ so that $\inf_{x\in S^1} d(\eta_k, \alpha\circ T_x) = \frac{\vn}{2}$.
Taking a subsequential limit of the $\eta_k$ yields a contradiction. 
  
\noindent {\em Proof:}

Suppose for the sake of contradiction that $\beta^t$ does not converge to $\alpha$. Then by Proposition \ref{DONE}, there is $\delta>0$ and a sequence of times $t'_j \ra \infty$
so that $\sup_{s\in S^1}d(\beta^{t'_j}(s), \alpha(S^1)) > \delta.$ Passing to a subsequence, we may assume $\beta^{t'_j}$ converges uniformly to a constant speed piecewise geodesic loop
$\gamma$ as $j\ra \infty$.

Since $\sup_{s\in S^1}d(\beta^{t'_j}(s), \alpha(S^1))>\delta$
for all $j$, $\sup_{s\in S^1}d(\gamma(s), \alpha(S^1))\geq\delta,$ so $\gamma$ is not a reparametrization of $\alpha$. But since $\beta^{t_j}\ra \alpha$,
$\lim_{t\ra \infty} l(t)=L$, and thus $\gamma$ has length $L$. So by hypothesis, $d(\gamma\circ T_x,\alpha)\geq \vn$ for all $x\in S^1$,
which is to say that $\inf_{x\in S^1} d(\gamma\circ T_x, \alpha)\geq \vn$.
We may assume that $t_j < t'_j$, $d(\beta^{t_j}, \alpha) < \frac{\vn}{2},$ and
$d(\beta^{t'_j}, \gamma) < \frac{\vn}{2}$ for all $j$. Then 
$$\inf_{x\in S^1} d(\beta^{t_j}\circ T_x, \alpha) \leq d(\beta^{t_j}, \alpha) < \frac{\vn}{2},$$ and 
$$\inf_{x\in S^1} d(\beta^{t'_j}\circ T_x, \alpha)\geq \lp \inf_{x\in S^1}  d(\gamma\circ T_x, \alpha)\rp - d(\beta^{t'_j}, \gamma) > \frac{\vn}{2}.$$
Since $(t,x)\mapsto d(\beta^t\circ T_x, \alpha)$ is continuous and $S^1$ is compact, $t\mapsto \inf_{x\in S^1} d(\beta^{t}\circ T_x, \alpha)$ is continuous.
Thus, for each $j$, there is $t''_j \in (t_j, t'_j)$,
so that $\inf_{x\in S^1} d(\beta^{t''_j}\circ T_x, \alpha) = \frac{\vn}{2}$. A subsequence of $\beta^{t''_j}$ converges to a constant speed loop $\eta$ of length $L$
with $$\inf_{x\in S^1} d(\eta\circ T_x, \alpha) = \frac{\vn}{2}.$$ Hence, $\eta$ is not a reparametrization of $\alpha$,
yet there is an $x$ so that $d(\eta\circ T_x, \alpha) = \frac{\vn}{2}<\vn$. This is a contradiction. \qed

Notation as above we have,
\begin{Corollary}\label{CURV}
Let $\alpha$ be a closed geodesic of length $L$ such that all sectional curvatures are negative at each point of the image of $\alpha$.
If $\beta^{t_j}\ra \alpha$ for some sequence $t_j\ra \infty$, then $\beta^t\ra \alpha$. 
\end{Corollary}
\noindent {\em Proof:} From the formula for second variation of arc-length, we know that $\alpha$ is isolated in the space of loops of length $L$, i.e. that $\alpha$
satisfies the hypotheses of Proposition \ref{ISO}.
\qed

{\em Remark:}
Suppose a subsequence $\beta^{t_j}$ converges to a closed geodesic $\alpha$.
Proposition \ref{ISO} shows that if $\alpha$ is isolated in the space of {\em rectifiable loops} 
of its length, then $\beta^t\ra \alpha$. Suppose instead $\alpha$ is merely isolated in the space of {\em closed geodesics} of its length. Suppose in addition that 
$\theta_i(t)$ converges to 0 for all $i$. Then all subsequential limits of $\beta^t$ are geodesics, so arguing as in the proof of Proposition \ref{ISO}, we can show $\beta^t \ra \alpha$. 
In particular, if the following conjecture holds, then $\beta^t$ converges to a closed geodesic (or a point) for pursuit on any compact manifold whose space of closed geodesics is discrete.
\begin{Conjecture}\label{ANGLE}
Let $M$ be a compact manifold, $\{b_i(t)\}_{i\in \Z/n}$ pursuit curves on $M$, $l_i(t)$ the associated lengths, $\theta_i(t)$ the associated angles. 
If $l_i(t)>0$ for all $i\in \Z/n$ and $t\geq 0$, then $\theta_i(t)\ra 0$.
\end{Conjecture}

Corollary \ref{CURV} and Proposition \ref{SS} imply Conjecture \ref{C} for compact manifolds of negative curvature.
Next, we prove Conjecture \ref{C} for manifolds of nonpositive curvature. First, a lemma:

\begin{Lemma}\label{NC}
Let $\alpha$ be a closed geodesic on a Riemannian manifold $M$. Suppose $\alpha(S^1)$ is contained in a non-positively curved open submanifold $U\subset M$. Suppose further that there is an $r>0$ so that $\inj(p)\geq r$ for all $p\in U$. Fix $\vn>0$. For sufficiently small $\delta$, the following property holds: for any 
$p_1,p_2 \in \ol{N}_{\delta}(\alpha)$
with $d(p_1,p_2) < r-\vn$, we have $[p_1,p_2]\subset \ol{N}_{\delta}(\alpha)$. 
\end{Lemma}

\noindent {\em Proof:}

We will show that it suffices to take $\delta$ small enough that 
\begin{enumerate}[(i)]
\item \label{LESS}
$\delta<\frac{\vn}{4}$
\item \label{IN}
for any $p,q \in \ol{N}_{\delta}(\alpha)$
with $d(p,q) < r-\frac{\vn}{4}$, we have $[p,q]\subset U$.
\end{enumerate}
(For \eqref{IN}, we use the continuous dependence of $[p,q]$ on $p,q$.)

Now, take $p_1,p_2 \in \ol{N}_{\delta}(\alpha)$ with $d(p_1,p_2) < r-\vn$. Let $\alpha:[0,1] \ra M$ be the shortest geodesic connecting $p_1$ to $p_2$.  
Choose $q_i$ in the image of $\alpha$ with $d(p_i, q_i)\leq \delta$ for $i=1,2$.
By condition \eqref{LESS} and the triangle inequality, $d(q_1,q_2) < r-\frac{\vn}{2}$,
so there is a unique shortest geodesic $\gamma:[0,1]\ra M$ connecting $q_1$ to $q_2$.

Observe that $d(\alpha(t), \gamma(t)) < r-\frac{\vn}{4}$ for all $t\in [0,1]$. Indeed, for $t\in [0,1/2]$, the path consisting of segments $[\alpha(t), p_1], [p_1, q_1], [q_1, \gamma(t)]$
has length less than $r-\frac{\vn}{4}$. Similarly, for $t\in [1/2, 1]$, the path consisting of $[\alpha(t), p_2], [p_2, q_2], [q_2, \gamma(t)]$ has length less than $r-\frac{\vn}{4}$. 

By \eqref{IN}, it follows that $[\alpha(t), \gamma(t)] \subset U$ for $t\in [0,1]$. The fact that $d(\alpha(t), \gamma(t)) < r-\frac{\vn}{4}$ for all $t\in [0,1]$ also implies
that the geodesic $[\alpha(t), \gamma(t)]$ varies smoothly in $t$. Since $U$ is nonpositively curved, 
we may apply the formula for the second variation of energy to the family of geodesics
$[\alpha(t), \gamma(t)]$ to conclude that $d^2(\alpha(t), \gamma(t))$ is convex as a function of $t$.
Therefore, we have for $t\in [0,1]$ that
$$d(\alpha(t), \gamma(t)) \leq \max\lbrac d(\alpha(0), \gamma(0)), d(\alpha(1), \gamma(1)) \rbrac = \max\lbrac d(p_1, q_1), d(p_2, q_2) \rbrac \leq \delta.$$
Thus, $[p_1,p_2] \subset \ol{N}_{\delta}(\alpha)$. \qed 

\begin{Proposition}\label{SUP}
If pursuit on a nonpositively curved compact manifold $M$ does not end in finite time, there is a closed geodesic $\alpha$ so
that $\sup_{s\in S^1}d(\beta^t(s), \alpha(S^1))\ra 0$ as $t\ra\infty$.
\end{Proposition}

\noindent {\em Proof:}

By Proposition \ref{SS}, there is a closed geodesic $\alpha$ and a sequence $t_j \ra \infty$ so that 
$\beta^{t_j} \ra \alpha$ uniformly.

Let $\vn>0$ be such that $d(b_i(0), b_{i+1}(0)) < \inj(M) - \vn$ for all $i$.
Take $\delta < \frac{\vn}{4}$. Then by the proof of Lemma \ref{NC}, if $p_1,p_2 \in \ol{N}_\delta(\alpha)$ and $d(p_1,p_2)<\inj(M)-\vn$, we have
$[p_1,p_2]\subset \ol{N}_\delta(\alpha)$. Also take $\delta$ small enough that $\ol{N}_\delta(\alpha)$ is a closed manifold with boundary, smoothly embedded in $M$.
Since $\beta^{t_j} \ra \alpha$, we have that the $b_i(t_J)$ are all in $\ol{N}_\delta(\alpha)$ for some sufficiently large $J$. Now, by Proposition \ref{CO},
we have $b_i(t) \in \ol{N}_\delta(\alpha)$ for all $i\in \Z/n$, $t\geq t_J$. By Lemma \ref{NC}, $\beta^t \subset \ol{N}_\delta(\alpha)$ for $t\geq t_J$.
Thus, $\sup_{s\in S^1}d(\beta^t(s), \alpha) \ra 0$ as $t\ra \infty$. \qed

As a consequence of Proposition \ref{DONE} and Proposition \ref{SUP}, we have Conjecture \ref{C} for manifolds of nonpositive curvature:
\begin{Theorem}
Let $M$ be a compact manifold of nonpositive sectional curvature. Suppose pursuit on $M$ with initial positions $\{b_i(0)\}_{i\in\Z/n}$ does not end in finite time,
and let $\beta^t$ be the associated family of piecewise geodesic loops.
Then there is a closed geodesic $\alpha$ so that $\beta^t \ra \alpha$.
\end{Theorem}

As an improvement of Corollary \ref{CURV}, we have the following result, which states that if $\beta^t$ gets close enough to a geodesic along which all sectional curvatures
are negative, then $\beta^t$ converges to that geodesic:
\begin{Proposition}\label{DUMBBELL}
Let $M$ be a Riemannian manifold with $\inj(M)>0$, and let $\alpha$ be a closed geodesic such that all sectional curvatures are negative at each point of the image of $\alpha$.
Fix $\vn>0$. Then there is a $\delta>0$ so that, if $d(b_i(0), b_{i+1}(0)) < \inj(M)-\vn$ for all $i$ and $\beta^{t_0}$ is uniformly $\delta$-close to $\alpha\circ T_x$ for some $t_0>0$ and some $x\in \R/\Z$, then
$\beta^t \ra \alpha$ as $t\ra \infty$. 
\end{Proposition}

\noindent {\em Proof:}

Take $\delta$ small enough so that
\be[(i)]
\item \label{ONE}
 if $p_1,p_2 \in \ol{N}_{\delta}(\alpha)$  and $d(p_1,p_2) < \inj(M)-\vn$, then
$[p_1,p_2]\subset \ol{N}_{\delta}(\alpha)$.  
\item \label{TWO}
 $\ol{N}_\delta(\alpha)$ is a closed manifold with boundary, smoothly embedded in $M$.
 \item \label{THREE}
 any loop uniformly $\delta$-close to $\alpha\circ T_x$ for some $x\in \R/\Z$ is homotopic to $\alpha$ through a family of loops in $\ol{N}_\delta(\alpha)$. 
\item \label{FOUR}
 any closed geodesic $\gamma$ in $\ol{N}_{\delta}(\alpha)$ homotopic to $\alpha$ through a family of loops in $\ol{N}_\delta(\alpha)$ differs from $\alpha$ by a rotation in the domain.
\ee
For condition \eqref{FOUR}, we argue on general grounds that taking $\delta$ small forces $\gamma$ to be uniformly close to $\alpha\circ T_x$ for some $x\in S^1$, and then we use the fact that a closed geodesic on $M$ along which all sectional curvatures are negative is isolated in the space of closed geodesics on $M$; the argument is straightforward, and we omit the details.
For \eqref{ONE}, we use Lemma \ref{NC}. For \eqref{TWO} and \eqref{THREE}, we take the image under $\exp$
of a neighborhood of the zero section in the normal bundle of $\alpha(S^1)$.

Now, suppose we have initial conditions for pursuit $\{b_i(0)\}_{i\in \Z/n}$ with $d(b_i(0), b_{i+1}(0)) < \inj(M)-\vn$ for all $i$, and the associated piecewise geodesic loop
$\beta^{t_0}$ is uniformly $\delta$-close to $\alpha\circ T_x$ for some $t_0>0$, $x\in \R/\Z$. By \eqref{ONE}, \eqref{TWO}, and Proposition \ref{CO}, we have $\beta^t \subset \ol{N}_\delta(\alpha)$
for all $t\geq t_0$. Using Proposition \ref{SS}, we get a sequence $t_j\ra \infty$ and a geodesic $\gamma$ contained in $\ol{N}_\delta(\alpha)$ so that $\beta^{t_j} \ra \gamma$ uniformly.
So $\beta^{t_0}$ is homotopic through a family of loops in $\ol{N}_\delta(\alpha)$ to $\gamma$. But by \eqref{THREE}, $\beta^{t_0}$ is also homotopic through a family of loops in $\ol{N}_\delta(\alpha)$ to $\alpha$.
Now by \eqref{FOUR}, $\gamma$ differs from $\alpha$ by a rotation in the domain. So by Corollary \ref{CURV}, $\beta^t \ra \alpha$ as $t\ra \infty$. 

\qed

\appendix
\section{Appendix}
We prove a result (Corollary \ref{AN} below) needed for the proof of Proposition \ref{FT2}.

\begin{Proposition}\label{AP}
Let $(M^n,g)$ be a Riemannian manifold, $p$ a point in $M$, $(U, x^i)$ a normal coordinate neighborhood centered at $p$.
Let $\|\cdot\|$ be the Euclidean norm on $(U, x^i)$.
Then for every $\vn>0$, there is an $r$ such that $B_r(p)\subset U$ 
and for every geodesic $\gamma:[0,a]\ra B_r(p)$, $\|\frac{\dot{\gamma}(0)}{\|\dot{\gamma}(0)\|} - \frac{\gamma(t) - \gamma(0)} {\|\gamma(t)-\gamma(0) \|} \|<\vn$,
for all $t\in (0,a]$.
  
\end{Proposition}

\noindent {\em Proof:}

Fix $\vn>0$.
Let $V$ be an open neighborhood of $p$ with closure contained in $U$. Then the  
Christoffel symbols $\Gamma_{ij}^k$ associated to $(U, x^i)$ are bounded on $V$.
Find $\mu$ so that $\lpi \Gamma_{ij}^k\rpi < \mu$ on $V$ for all $i,j,k$. Take $r$ small enough that $B_r(p) \subset V$.
Since $g_{ij}(p) = \delta_{ij}$, we may take $r$
small enough that any vector in $TB_r(p)$ of unit length with respect to $g$ has length less than 2
with respect to the Euclidean norm. Now, if $\gamma:[0,a]\ra B_r(p)$ is a unit speed geodesic, we have
$$\frac{d^2\gamma^k}{dt^2} = -\Gamma^{k}_{ij}\frac{d\gamma^i}{dt}\frac{d\gamma^j}{dt},$$
so $\lpi \frac{d^2\gamma^k}{dt^2} \rpi < 4\mu n^2$ and thus $\|\frac{d^2\gamma}{dt^2}\| < 4\mu n^{\frac{5}{2}}$.
Integrating, we find for each $t\in [0,a]$ that
$$\| \dot{\gamma}(t)-\dot{\gamma}(0)\|<4t\mu n^{\frac{5}{2}} \leq 8r\mu n^{\frac{5}{2}},$$
so taking $r$ is smaller than $\lbrac 8\mu n^{\frac{5}{2}}\rbrac^{-1}\vn$, we have 
\begin{equation}\label{DIR}
\|\dot{\gamma}(t)-\dot{\gamma}(0)\|<\vn.
\end{equation}
Integrating again, we get
$\|\gamma(t) - \gamma(0) - t\dot{\gamma}(0)\|<t\vn$ and so $ \| \frac{\gamma(t)-\gamma(0)}{t}-\dot{\gamma}(0)\| < \vn$ for $t\in (0,a]$.
Assuming $r$ is chosen small enough so that any $v\in TB_r(p)$ with $g(v,v)=1$ satisfies $\| v- \frac{v}{\|v\|}\|<\vn$,
we have  $\left\| \frac{\gamma(t)-\gamma(0)}{t}-\frac{\dot{\gamma}(0)}{\|\dot{\gamma}(0)\|}\right\| < 2\vn$.
Hence, the Euclidean distance from the unit vector $\frac{\dot{\gamma}(0)}{\|\dot{\gamma}(0)\|}$ to the line $\R \lp \gamma(t)-\gamma(0)\rp$ is less than $2\vn$.
For $\vn$ small enough, this implies
$$  \left \| \frac{\gamma(t)-\gamma(0)}{\| \gamma(t)-\gamma(0) \|} - \frac{\dot{\gamma}(0)}{\|\dot{\gamma}(0)\|} \right\|< 3\vn .$$ 
\qed

Notation as in the last Proposition, we have the following
\begin{Corollary}\label{AN}
For every $\vn>0$, there is an $r$ so that $B_r(p)\subset U$ and for any two geodesics
$\gamma_1:[0,a_1] \ra B_r(p),~\gamma_2:[0, a_2] \ra B_r(p)$ with $\gamma_1(0) = \gamma_2(0)$,
the metric angle between $\dot{\gamma}_1(0), \dot{\gamma}_2(0)$
is within $\vn$ of the Euclidean angle between $\gamma_1(a_1) - \gamma_1(0)$ and $\gamma_2(a_2) - \gamma_2(0)$.  
\end{Corollary}

\noindent {\em Proof:}

By the last part and uniform continuity of the spherical distance function $S^{n-1}\times S^{n-1}\ra \R$ on the unit sphere,
we can choose $r$ small enough that the Euclidean angle between $\gamma_1(a_1) - \gamma_1(0)$ and $\gamma_2(a_2) - \gamma_2(0)$
is within $\frac{\vn}{2}$ of the Euclidean angle between $\dot{\gamma}_1(0), \dot{\gamma}_2(0)$, for any two unit speed geodesics
$\gamma_1:[0,a_1] \ra B_r(p),~\gamma_2:[0, a_2] \ra B_r(p)$ with $\gamma_1(0) = \gamma_2(0)$. Then, if necessary, we choose $r$ smaller so that
for any two vectors $u,v\in TB_r(p)$ based at the same point, the Euclidean angle between $u,v$ is within $\frac{\vn}{2}$ of the metric angle. \qed

\section*{Acknowledgements}
This research was conducted mostly at the SUMMER@ICERM Undergraduate Summer Research Program in 2012.
I would like to thank Tarik Aougab and Sergei Tabachnikov for their mentorship.
I would like to thank Francisc Bozgan for pointing out the application of Jensen's inequality in Section \ref{EUC}. 
I would like to thank Anton Petrunin for a MathOverflow answer which helped with the proof of Proposition \ref{ISO}.

\end{document}